\documentclass[twoside,11pt]{article}

\usepackage{amsthm,amssymb,amsmath}
\usepackage{graphicx}

\newtheorem{thm}{Theorem}
\newtheorem{cor}{Corollary}

\theoremstyle{definition}
\newtheorem{defn}{Definition}
\theoremstyle{remark}
\newtheorem{rem}{Remark}
\numberwithin{equation}{section}
\begin{document}

\pagestyle{myheadings} \markboth{ \rm \centerline {R. N. Mohapatra and B. Szal}} {\rm \centerline { }}

\begin{titlepage}
\title{\bf {Degree of Convergence of an Integral Operator }}
\author {\bf\Large R. N. Mohapatra $^{1}$ and 
B. Szal $^{2}$\\
{\small $^1$ University of Central Florida}\\
{\small Orlando, FL 32816, USA}\\
{\small ramm@pegasus.cc.ucf.edu}\\
{\small $^2$ University of Zielona G\'{o}ra}\\
{\small Faculty of Mathematics, Computer Science and Econometrics}\\
{\small 65-516 Zielona G\'{o}ra, ul. Szafrana 4a, Poland}\\
{\small B.Szal@wmie.uz.zgora.pl} }
\end{titlepage}

\date{}
\maketitle
\begin{abstract}
In this paper we define an integral operator on $L^{p}$ and obtain its degree of convergence in the appropriate norm. By specializing the kernel of the integral operator we obtain many known results as corollaries. We have also applied our results to obtain results on singular integral operators.
\end{abstract}

\noindent{\it Keywords and phrases:} 

\noindent { \it 2000 Mathematics Subject Classification:}

\maketitle

\section{Introduction}

Let $L^{p}\equiv L^{p}\left( 
\mathbb{R}
\right) $ with fixed $1\leq p\leq \infty $ be the space of all real-valued
functions Lebesgue integrable to the $p-$th power over $%
\mathbb{R}
$ if $1\leq p<\infty $ and uniformly continuous and bounded on $%
\mathbb{R}
$ if $p=\infty $. We define the norm in $L^{p}$, as usual, by the formula%
\begin{equation}
\left\Vert f\right\Vert _{p}\equiv \left\Vert f\left( \cdot \right)
\right\Vert _{p}:=\left\{ 
\begin{array}{c}
\left\{ \int\limits_{%
\mathbb{R}
}\left\vert f\left( x\right) \right\vert ^{p}dx\right\} ^{\frac{1}{p}}\text{
\ \ if \ \ }1\leq p<\infty , \\ 
\underset{x\in 
\mathbb{R}
}{\sup }\left\vert f\left( x\right) \right\vert \text{ \ \ \ \ \ \ \ \ \ \ \
\ \ \ \ \ \ if \ \ }p=\infty \text{.}%
\end{array}%
\right.  \label{i00}
\end{equation}%
Denote by $\omega \left( f;\cdot \right) _{p}$ the modulus of continuity of $%
f\in L^{p}$, i.e.,%
\begin{equation*}
\omega \left( f;t\right) _{p}:=\underset{0\leq h\leq t}{\sup }\left\Vert
\Delta _{h}f\left( \cdot \right) \right\Vert _{p},\text{ \ \ }t\geq 0,
\end{equation*}%
where $\Delta _{h}f\left( x\right) =f\left( x+h\right) -f\left( x\right) $.

Consider the family of integral operators%
\begin{equation}
F_{\lambda }\left( f;x\right) :=\lambda \int\limits_{%
\mathbb{R}
}f\left( t\right) \mathcal{K} \left( \lambda \left( t-x\right) \right) dt%
\text{, \ \ }\lambda >0,  \label{i1}
\end{equation}%
with Fej\'{e}r type kernel $\mathcal{K} $ \cite[p. 126]{1}:%
\begin{equation}
\mathcal{K} \left( -x\right) =\mathcal{K} \left( x\right) ,  \label{i2}
\end{equation}%
\begin{equation}
\int\limits_{%
\mathbb{R}
}\mathcal{K} \left( x\right) dx=1,  \label{i3}
\end{equation}%
\begin{equation}
\underset{-1\leq x\leq 1}{\sup }\left\vert \mathcal{K} \left( x\right)
\right\vert <\infty ,  \label{i4}
\end{equation}%
\begin{equation}
\underset{x\in 
\mathbb{R}
}{\sup }~x^{2}\left\vert \mathcal{K} \left( x\right) \right\vert <\infty .
\label{i5}
\end{equation}

Under these conditions, the integral (\ref{i1}) represents a linear operator
acting from $L^{p}$ to $L^{p}.$

For fixed $m\in 
\mathbb{N}
\cup \left\{ 0\right\} $ and $1\leq p\leq \infty $, we denote by $L_{m}^{p}$
the set of all $f\in L^{p}$ whose derivatives $f^{\prime },f^{\prime \prime
},...,f^{\left( m\right) }$ also belong to $L^{p}$. The norm in these $%
L_{m}^{p}$ is defined by (\ref{i00}), i.e., for $f\in L_{m}^{p}$, we have $%
\left\Vert f\right\Vert _{p,~m}=\left\Vert f\right\Vert _{p}$, where $%
\left\Vert f\right\Vert _{p}$ is defined by (\ref{i00}). It is clear that $%
L_{0}^{p}\equiv L^{p}$.

\begin{defn}
Let $f\in L_{m}^{p}$ for fixed $m\in 
\mathbb{N}
\cup \left\{ 0\right\} $ and $1\leq p\leq \infty $. We define a family of
modified integral operators by the formula%
\begin{equation}
F_{\lambda ,~m}\left( f;x\right) :=\lambda \int\limits_{%
\mathbb{R}
}\sum\limits_{j=0}^{m}\frac{f^{\left( j\right) }\left( t\right) }{j!}\left(
x-t\right) ^{j}\mathcal{K} \left( \lambda \left( t-x\right) \right) dt
\label{i6}
\end{equation}%
for $x\in 
\mathbb{R}
$ and $\lambda >0$.

In particular, we have $F_{\lambda ,~0}\left( f;\cdot \right) \equiv
F_{\lambda }\left( f;\cdot \right) $ for $f\in L^{p}$.
\end{defn}

It is obvious that the formula (\ref{i6}) can be rewritten in the following form:%
\begin{equation*}
F_{\lambda ,~m}\left( f;x\right) :=\sum\limits_{j=0}^{m}\frac{\left(
-1\right) ^{j}}{j!}\lambda \int\limits_{%
\mathbb{R}
}f^{\left( j\right) }\left( t+x\right) t^{j}\mathcal{K} \left( \lambda
t\right) dt
\end{equation*}%
for every $f\in L_{m}^{p}$, $x\in 
\mathbb{R}
$ and $\lambda >0$.

If (\ref{i2}) holds and for any $j=0,1,2,...,m$%
\begin{equation*}
\int\limits_{0}^{\infty }u^{j}\left\vert \mathcal{K} \left( u\right)
\right\vert <\infty ,
\end{equation*}%
then for fixed $m\in 
\mathbb{N}
\cup \left\{ 0\right\} $ and $\lambda >0$ the integral (\ref{i6}) is a linear
operator from space $L_{m}^{p}$ into $L^{p}$ (see Remark 2).

Denote by $H^{\omega ^{\ast },~p}$ the set of all functions $f\in L^{p}\left(
1\leq p\leq \infty \right) $ satisfying the condition%
\begin{equation*}
\underset{h\neq 0}{\sup }~h^{\omega ^{\ast },~p}\left( f;h\right) <\infty ,
\end{equation*}%
where%
\begin{equation*}
h^{\omega ^{\ast },~p}\left( f;h\right) :=\frac{\left\Vert \Delta _{h}f\left(
\cdot \right) \right\Vert _{p}}{\omega ^{\ast }\left( \left\vert
h\right\vert \right) },\text{ \ \ \ \ \ }h^{0,~p}\left( f;h\right) =0
\end{equation*}%
and, for $t\geq 0$, $\omega ^{\ast }\left( t\right) $ is a nondecreasing
function. We can show that $H^{\omega ^{\ast },~p}$ is a Banach space with
respect to the generalized H\"{o}lder norm%
\begin{equation}
\left\Vert f\right\Vert _{\omega ^{\ast },~p}:=\left\Vert f\right\Vert _{p}+%
\underset{h\neq 0}{\sup }~ h^{\omega ^{\ast },~p}\left( f;h\right) .  \label{i0}
\end{equation}%
Suppose that $H^{\omega ,~p}$ is the set of functions $f\in L^{p}\left( 1\leq
p\leq \infty \right) $ satisfying the condition%
\begin{equation*}
\underset{h\neq 0}{\sup }~h^{\omega ,~p}\left( f;h\right) =\underset{h\neq 0}{%
\sup }\frac{\left\Vert \Delta _{h}f\left( \cdot \right) \right\Vert _{p}}{%
\omega \left( \left\vert h\right\vert \right) }<\infty
\end{equation*}%
and contained in the space $H^{\omega ^{\ast },~p}$, $H^{\omega ,~p}\subset
H^{\omega ^{\ast },~p}$, where, for $t\geq 0$, $\omega \left( t\right) $ is a
nondecreasing function. In particular, setting%
\begin{equation*}
\omega \left( t\right) =t^{\alpha }\text{, \ \ \ }\omega ^{\ast }\left(
t\right) =t^{\beta }\text{, \ \ \ }t\geq 0\text{ \ \ and \ \ }0\leq \beta
<\alpha \leq 1~,
\end{equation*}%
for $H^{\omega ^{\ast },~p}$ we obtain the spaces%
\begin{equation*}
H^{\beta ,~p}:=\left\{ f\in L^{p}:\omega \left( f,t\right) _{p}\leq
C_{1}~t^{\beta }\right\}
\end{equation*}%
with H\"{o}lder norm%
\begin{equation*}
\left\Vert f\right\Vert _{\beta ,~p}:=\left\Vert f\right\Vert _{p}+\underset{%
h\neq 0}{\sup }\frac{\left\Vert \Delta _{h}f\left( \cdot \right) \right\Vert
_{p}}{\left\vert h\right\vert ^{\beta }},
\end{equation*}%
and for the set $H^{\omega ,~p}$ we have%
\begin{equation*}
H^{\alpha ,~p}:=\left\{ f\in L^{p}:\omega \left( f,t\right) _{p}\leq
C_{2}~t^{\alpha }\right\} ,
\end{equation*}%
\begin{equation*}
H^{\alpha ,~p}\subset H^{\beta,~p}.
\end{equation*}

For fixed $m\in 
\mathbb{N}
\cup \left\{ 0\right\} $ and $1\leq p\leq \infty $, we denote by $%
H_{m}^{\omega ^{\ast },~p}$ (or $H_{m}^{\omega ,~p}$) the set of all $f\in
H^{\omega ^{\ast },~p}$ ($f\in H^{\omega ,~p}$) whose derivatives $f^{\prime
},f^{\prime \prime },...,f^{\left( m\right) }$ also belong to $H^{\omega
^{\ast },~p}$ (or $H^{\omega ,~p}$), where, for $t\geq 0$, $\omega ^{\ast
}\left( t\right) $ (or $\omega \left( t\right) $) is a nondecreasing
function. The norm in these $H_{m}^{\omega ^{\ast },~p}$ is defined by (\ref%
{i0}), i.e., for $f\in H_{m}^{\omega ^{\ast },~p}$, we have $\left\Vert
f\right\Vert _{\omega ^{\ast },~p,~m}=\left\Vert f\right\Vert _{\omega ^{\ast
},~p}$, where $\left\Vert f\right\Vert _{\omega ^{\ast },~p}$ is defined by (%
\ref{i0}). It is clear that $H_{0}^{\omega ^{\ast },~p}\equiv H^{\omega
^{\ast },~p}$.

Throughout the paper we shall use the following notation:%
\begin{equation*}
E_{\lambda }\left( x\right) =E_{\lambda }\left( f;x\right) =F_{\lambda
}\left( f;x\right) -f\left( x\right) ,
\end{equation*}%
\begin{equation*}
E_{\lambda }\left( x+h,x\right) =E_{\lambda }\left( f;x+h,x\right)
=E_{\lambda }\left( x+h\right) -E_{\lambda }\left( x\right)
\end{equation*}%
and%
\begin{equation*}
\phi _{x}\left( t\right) =f\left( x+t\right) +f\left( x-t\right) -2f\left(
x\right) .
\end{equation*}
The object of this paper is to study degree of convergence of the integral operator $F_{\lambda ,~m}\left( f\right)$ to $f$ in the appropriate norm and to deduce many interesting results as corollaries. We also apply our results to obtain degree of convergence of singular integrals. 

\section{Statement of the results}

\begin{thm}
Suppose that $0\leq \beta <\eta \leq 1$ and (\ref{i2})-(\ref{i5}) hold.
Then, for any $f\in H^{\omega ,~p}\left( 1\leq p\leq \infty \right) $, the
following relation is true%
\begin{equation}
\left\Vert F_{\lambda }\left( f\right) -f\right\Vert _{\omega ^{\ast
},~p}=O\left( 1\right) \underset{h\neq 0}{\sup }~\frac{\left( \omega \left(
\left\vert h\right\vert \right) \right) ^{\frac{\beta }{\eta }}}{\omega
^{\ast }\left( \left\vert h\right\vert \right) }\left\{ \frac{1}{\lambda }%
\int\limits_{1}^{\lambda }\left( \omega \left( \frac{1}{u}\right) \right)
^{1-\frac{\beta }{\eta }}du\right\} ,  \label{t1}
\end{equation}%
where $\lambda >1$.
\end{thm}

\begin{rem}
For $p=\infty $ (\ref{t1}) was proved in \cite[Theorem 1]{5}.
\end{rem}

\begin{thm}
Suppose that $0\leq \beta <\eta \leq 1,$ $\lambda \geq \lambda _{0}>0$ and (%
\ref{i2})-(\ref{i4}) hold. If 
\begin{equation}
\int\limits_{1}^{\infty }u\left\vert \mathcal{K} \left( u\right)
\right\vert du<\infty ,  \label{t2}
\end{equation}%
then for $f\in H^{\omega ,~p}\left( 1\leq p\leq \infty \right) $%
\begin{equation*}
\left\Vert F_{\lambda }\left( f\right) -f\right\Vert _{\omega ^{\ast
},~p}=O\left( 1\right) \underset{h\neq 0}{\sup }~\frac{\left( \omega \left(
\left\vert h\right\vert \right) \right) ^{\frac{\beta }{\eta }}}{\omega
^{\ast }\left( \left\vert h\right\vert \right) }\left( \omega \left( \frac{1%
}{\lambda }\right) \right) ^{1-\frac{\beta }{\eta }}.
\end{equation*}
\end{thm}

\begin{rem}
Let $m\in 
\mathbb{N}
\cup \left\{ 0\right\} $, $1\leq p\leq \infty $ and (\ref{i2}) holds. If for
any $j=0,1,2,...,m$%
\begin{equation}
\int\limits_{0}^{\infty }u^{j}\left\vert \mathcal{K} \left( u\right)
\right\vert du<\infty ,  \label{t3}
\end{equation}%
then for every $f\in L_{m}^{p}$ and $\lambda >0$ we have%
\begin{equation*}
\left\Vert F_{\lambda ,~m}\left( f\right) \right\Vert _{p}=O\left( 1\right)
\sum\limits_{j=0}^{m}\frac{\left\Vert f^{\left( j\right) }\right\Vert _{p}}{%
j!\lambda ^{j}}.
\end{equation*}%
Moreover, if $f\in H_{m}^{\omega ^{\ast },~p}$ and $\lambda >0$, then%
\begin{equation*}
\left\Vert F_{\lambda ,~m}\left( f\right) \right\Vert _{\omega ^{\ast
},~p}=O\left( 1\right) \sum\limits_{j=0}^{m}\frac{\left\Vert f^{\left(
j\right) }\right\Vert _{\omega ^{\ast },~p}}{j!\lambda ^{j}}.
\end{equation*}
\end{rem}

\begin{thm}
Suppose that $0\leq \beta <\eta \leq 1$, $m\in 
\mathbb{N}
$, $\lambda \geq \lambda _{0}>0$ and (\ref{i2})-(\ref{i3}) hold. If for any $%
j=1,2,...,m+1$%
\begin{equation}
\int\limits_{0}^{\infty }u^{j}\left\vert \mathcal{K} \left( u\right)
\right\vert du<\infty ,  \label{t4}
\end{equation}%
then for $f\in H_{m}^{\omega ,~p}\left( 1\leq p\leq \infty \right) $%
\begin{equation*}
\left\Vert F_{\lambda ,~m}\left( f\right) -f\right\Vert _{\omega ^{\ast
},~p}=O\left( 1\right) \underset{h\neq 0}{\sup }~\frac{\left( \omega \left(
\left\vert h\right\vert \right) \right) ^{\frac{\beta }{\eta }}}{\omega
^{\ast }\left( \left\vert h\right\vert \right) }\frac{1}{\lambda ^{m}}\left(
\omega \left( \frac{1}{\lambda }\right) \right) ^{1-\frac{\beta }{\eta }}.
\end{equation*}
\end{thm}

Setting%
\begin{equation*}
\omega \left( t\right) =t^{\alpha }\text{, \ \ }\omega ^{\ast }\left(
t\right) =t^{\beta }\text{, \ \ }0\leq \beta <\alpha \leq 1\text{, \ \ }\eta
=\alpha
\end{equation*}%
in the assumptions of Theorem 1, 2 and 3, we obtain the following
corollaries:

\begin{cor}
Suppose that $0\leq \beta <\alpha \leq 1$ and (\ref{i2})-(\ref{i5}) hold.
Then, for any $f\in H^{\alpha ,~p}\left( 1\leq p\leq \infty \right) $, the
following relation is true%
\begin{equation*}
\left\Vert F_{\lambda }\left( f\right) -f\right\Vert _{\beta ,~p}=\left\{ 
\begin{array}{c}
O\left( \lambda ^{\beta -\alpha }\right) \text{ \ \ if \ \ }\alpha -\beta <1,
\\ 
O\left( \frac{\ln \lambda }{\lambda }\right) \text{ \ \ if \ \ }\alpha
-\beta =1,%
\end{array}%
\right.
\end{equation*}%
where $\lambda >1$.
\end{cor}

\begin{cor}
Suppose that $0\leq \beta <\alpha \leq 1,$ $\lambda \geq \lambda _{0}>0$ and
(\ref{i2})-(\ref{i4}) hold. If 
\begin{equation*}
\int\limits_{1}^{\infty }u\left\vert \mathcal{K} \left( u\right)
\right\vert du<\infty ,
\end{equation*}%
then for $f\in H^{\alpha ,~p}\left( 1\leq p\leq \infty \right) $%
\begin{equation*}
\left\Vert F_{\lambda }\left( f\right) -f\right\Vert _{\beta ,~p}=O\left(
\lambda ^{\beta -\alpha }\right) .
\end{equation*}
\end{cor}

\begin{cor}
Suppose that $0\leq \beta <\alpha \leq 1$, $m\in 
\mathbb{N}
$, $\lambda \geq \lambda _{0}>0$ and (\ref{i2})-(\ref{i3}) hold. If for any $%
j=1,2,...,m+1$%
\begin{equation*}
\int\limits_{0}^{\infty }u^{j}\left\vert \mathcal{K} \left( u\right)
\right\vert du<\infty ,
\end{equation*}%
then for $f\in H_{m}^{\alpha ,~p}\left( 1\leq p\leq \infty \right) $%
\begin{equation*}
\left\Vert F_{\lambda ,~m}\left( f\right) -f\right\Vert _{\beta ,~p}=O\left(
\lambda ^{\beta -\alpha -m}\right) .
\end{equation*}
\end{cor}

\section{Examples}

\subsection{The Riesz means of the Fourier series}

Suppose that the kernel $\mathcal{K} $ satisfies the conditions (\ref{i2})-(%
\ref{i5}). Then, as is well known (see \cite[p. 132]{1}), if we consider $%
2\pi $ periodic function  $f$ $\in L^{p}$ with the Fourier series%
\begin{equation*}
S\left( f\right) =\sum\limits_{k=-\infty }^{\infty }c_{k}\left( f\right)
e^{ikx},
\end{equation*}%
the family of operators of Fej\'{e}r type (\ref{i1}) can be transformed into
sequences of linear means of series of the function $f$%
\begin{equation*}
F_{n}\left( f;x\right) =\sum\limits_{k=-\infty }^{\infty }\varphi \left( 
\frac{k}{n}\right) c_{k}\left( f\right) e^{ikx}\text{, \ \ }n\in 
\mathbb{N}
,
\end{equation*}%
where the values of $\varphi \left( \frac{k}{n}\right) $ coincide for $t=%
\frac{k}{n}$ with the values of the function $\varphi \left( t\right) $,
which is the Fourier transform of the kernel $\mathcal{K} $.

Consider the Riesz means of the Fourier series $S\left( f\right) :$ 
\begin{equation*}
R_{n}\left( \gamma ,f;x\right) =\sum_{k=-n}^{n}\left( 1-\frac{\left\vert
k\right\vert }{n}\right) ^{\gamma }c_{k}\left( f\right) e^{ikx}\text{, \ \ }%
\gamma >0.
\end{equation*}%
For this mean 
\begin{equation*}
\varphi _{\gamma }\left( t\right) =\left\{ 
\begin{array}{c}
\left( 1-\left\vert t\right\vert \right) ^{\gamma }\text{, \ }\left\vert
t\right\vert \leq 1, \\ 
0\text{, \ \ \ \ }\left\vert t\right\vert \geq 1,%
\end{array}%
\right. 
\end{equation*}%
and consequently%
\begin{equation*}
\mathcal{K} _{R\left( \gamma \right) }\left( t\right) =\frac{1}{2\pi }%
\int\limits_{%
\mathbb{R}
}\varphi _{\gamma }\left( x\right) e^{itx}dx
\end{equation*}%
\begin{equation}
=\frac{1}{\pi }\left\{ \cos t\sum\limits_{k=0}^{\infty }\frac{\left(
-1\right) ^{k}t^{2k}}{\left( 2k\right) !\left( 2k+1+\gamma \right) }+\sin
t\sum\limits_{k=0}^{\infty }\frac{\left( -1\right) ^{k}t^{2k+1}}{\left(
2k+1\right) !\left( 2k+2+\gamma \right) }\right\} .  \label{t5}
\end{equation}%
It is clear that the kernel $\mathcal{K} _{R\left( \gamma \right) }$
satisfies the conditions (\ref{i2}), (\ref{i4}) and (\ref{i5}). Moreover,
the function $\varphi _{\gamma }$ is the Fourier transform of the kernel $%
\mathcal{K} _{R\left( \gamma \right) }$, i.e.,%
\begin{equation*}
\varphi _{\gamma }\left( t\right) =\int\limits_{%
\mathbb{R}
}\mathcal{K} _{R\left( \gamma \right) }\left( x\right) e^{-itx}dx
\end{equation*}%
for all $t\in 
\mathbb{R}
$.  Thus%
\begin{equation*}
1=\varphi _{\gamma }\left( 0\right) =\int\limits_{%
\mathbb{R}
}\mathcal{K} _{R\left( \gamma \right) }\left( x\right) dx
\end{equation*}%
and the condition (\ref{i3}) is valid, too. Hence, by Theorem 1, we obtain
the following results.

\begin{cor}
Suppose that $0\leq \beta <\eta \leq 1$ and $\gamma >0$. Then, for any $2\pi $ periodic
function $f\in H^{\omega ,~p}\left( 1\leq p\leq \infty \right) $, the
following relation is true:%
\begin{equation*}
\left\Vert R_{n}\left( \gamma ,f\right) -f\right\Vert _{\omega ^{\ast
},~p}=O\left( 1\right) \underset{h\neq 0}{\sup }~\frac{\left( \omega \left(
\left\vert h\right\vert \right) \right) ^{\frac{\beta }{\eta }}}{\omega
^{\ast }\left( \left\vert h\right\vert \right) }\left\{ \frac{1}{n}%
\sum\limits_{k=1}^{n}\left( \omega \left( \frac{1}{k}\right) \right) ^{1-%
\frac{\beta }{\eta }}\right\} .
\end{equation*}
\end{cor}

\begin{cor}
Suppose that $0\leq \beta <\alpha \leq 1$ and $\gamma >0$. Then, for any $2\pi $ periodic
function $f\in H^{\alpha ,~p}\left( 1\leq p\leq \infty \right) $, the
following relation holds:%
\begin{equation*}
\left\Vert R_{n}\left( \gamma ,f\right) -f\right\Vert _{\beta ,~p}=\left\{ 
\begin{array}{c}
O\left( n^{\beta -\alpha }\right) \text{ \ \ if \ \ }\alpha -\beta <1, \\ 
O\left( \frac{\ln n}{n}\right) \text{ \ \ if \ \ }\alpha -\beta =1.%
\end{array}%
\right. 
\end{equation*}
\end{cor}

In particular, putting $\gamma =1$ in (\ref{t5}), we get the Fej\'{e}r
kernel 
\begin{equation*}
\mathcal{K} _{\sigma }\left( t\right) :=\mathcal{K} _{R\left( 1\right)
}\left( t\right) =\frac{2}{\pi }\left( \frac{\sin \left( t/2\right) }{t}%
\right) ^{2}
\end{equation*}%
and consequently the Fej\'{e}r means of the Fourier series $S\left( f\right) 
$ is given by 
\begin{equation*}
\sigma _{n}\left( f;x\right) :=R_{n}\left( 1,f;x\right)
=\sum\limits_{k=-\infty }^{\infty }\varphi \left( \frac{\kappa }{n}\right)
c_{k}\left( f\right) e^{ikx}=\sum\limits_{k=-n}^{n}\left( 1-\frac{%
\left\vert k\right\vert }{n}\right) c_{k}\left( f\right) e^{ikx}.
\end{equation*}%
Hence, from the above corollaries, we obtain the results from \cite{3}, \cite{5}
and \cite{9}.

\subsection{The Poisson operator}

Let%
\begin{equation*}
\mathcal{K} _{\overline{P}}\left( t\right) =\frac{1}{\pi }\frac{1}{1+t^{2}}%
\text{, \ \ }\lambda =\frac{1}{\varepsilon }\text{, \ \ }\varepsilon >0.
\end{equation*}%
Then we obtain from (\ref{i1}) the Poisson singular integral of a function $f
$ $\in L^{p}$, i.e%
\begin{equation*}
\overline{P}_{\varepsilon }\left( f;x\right) =\frac{\varepsilon }{\pi }%
\int\limits_{%
\mathbb{R}
}f\left( x+t\right) \frac{1}{\varepsilon ^{2}+t^{2}}dt,
\end{equation*}%
The approximation properties of this integral was given in \cite{11}%
, for example. It is clear that the kernel $\mathcal{K} _{\overline{P}}$
satisfies conditions (\ref{i2})-(\ref{i5}). Hence, by
Theorem 1, we obtain the following assertion.

\begin{cor}
Suppose that $0\leq \beta <\eta \leq 1$. Then, for any $2\pi $ periodic
function $f\in H^{\omega ,~p}\left( 1\leq p\leq \infty \right) $, the
following relation is true%
\begin{equation*}
\left\Vert \overline{P}_{\varepsilon }\left( f\right) -f\right\Vert _{\omega
^{\ast },~p}=O\left( 1\right) \underset{h\neq 0}{\sup }~\frac{\left( \omega
\left( \left\vert h\right\vert \right) \right) ^{\frac{\beta }{\eta }}}{%
\omega ^{\ast }\left( \left\vert h\right\vert \right) }\left\{ \varepsilon
\int\limits_{1}^{1/\varepsilon }\left( \omega \left( \frac{1}{u}\right)
\right) ^{1-\frac{\beta }{\eta }}du\right\} ,
\end{equation*}%
where $0<\varepsilon <1$.
\end{cor}

\begin{cor}
Suppose that $0\leq \beta <\alpha \leq 1$. Then, for any $2\pi $ periodic
function $f\in H^{\alpha ,~p}\left( 1\leq p\leq \infty \right) $, the
following relation is true%
\begin{equation*}
\left\Vert \overline{P}_{\varepsilon }\left( f\right) -f\right\Vert _{\beta
,~p}=\left\{ 
\begin{array}{c}
O\left( \varepsilon ^{\beta -\alpha }\right) \text{ \ \ if \ \ }\alpha
-\beta <1, \\ 
O\left( \varepsilon \ln \left( 1/\varepsilon \right) \right) \text{ \ \ if \
\ }\alpha -\beta =1,%
\end{array}%
\right. 
\end{equation*}%
as $\varepsilon \rightarrow 0^{+}$.
\end{cor}

\subsection{The Picard and Gauss-Weierstrass operators}

Taking 
\begin{equation*}
\mathcal{K} _{P}\left( t\right) =\frac{1}{2}\exp \left( -\left\vert
t\right\vert \right) \text{, \ \ }\lambda =\frac{1}{r},\text{ \ \ }r>0
\end{equation*}%
and%
\begin{equation*}
\mathcal{K} _{W}\left( t\right) =\frac{1}{\sqrt{\pi }}\exp \left(
-t^{2}\right) \text{, \ \ }\lambda =\frac{1}{\sqrt{2r}}\text{, \ \ }r>0
\end{equation*}%
we obtain from (\ref{i1}), respectively, the Picard singular integral and
the Gauss-Weierstrass singular integrals of a function $f$ $\in L^{p}$, i.e.,%
\begin{equation*}
P_{r}\left( f;x\right) =\frac{1}{2r}\int\limits_{%
\mathbb{R}
}f\left( x+t\right) \exp \left( \frac{-\left\vert t\right\vert }{r}\right)
dt,
\end{equation*}%
\begin{equation*}
W_{r}\left( f;x\right) =\frac{1}{\sqrt{4\pi r}}\int\limits_{%
\mathbb{R}
}f\left( x+t\right) \exp \left( -\frac{t^{2}}{4r}\right) dt.
\end{equation*}%
The limit properties (as $r\rightarrow 0^{+}$) of these integrals were given
in many papers and monographs (see, e.g., \cite{2,3, 4}).

It is clear (see \cite[Lemma 1]{5}) that for every $m\in 
\mathbb{N}
\cup \left\{ 0\right\} $%
\begin{equation*}
\int\limits_{0}^{\infty }u^{m}\left\vert \mathcal{K} _{P}\left( u\right)
\right\vert du=\frac{m!}{2}
\end{equation*}%
and%
\begin{equation*}
\int\limits_{0}^{\infty }u^{m}\left\vert \mathcal{K} _{W}\left( u\right)
\right\vert du=\left\{ 
\begin{array}{c}
\frac{1}{2} \\ 
\frac{\left( 2k-1\right) !!}{2^{k+1}} \\ 
\frac{k!}{2\sqrt{\pi }}%
\end{array}%
\begin{array}{c}
\text{if} \\ 
\text{if} \\ 
\text{if}%
\end{array}%
\begin{array}{c}
m=0, \\ 
m=2k\geq 2, \\ 
m=2k+1\geq 1.%
\end{array}%
\right.
\end{equation*}%
Hence, by Theorems 2 and 3, we obtain the following assertion.

\begin{cor}
Suppose that $0\leq \beta <\eta \leq 1$, $m\in 
\mathbb{N}
\cup \left\{ 0\right\} $ and $f\in H_{m}^{\omega ,~p}\left( 1\leq p\leq
\infty \right) $. Then%
\begin{equation*}
\left\Vert P_{r}\left( f\right) -f\right\Vert _{\omega ^{\ast },~p}=O\left(
1\right) \underset{h\neq 0}{\sup }\frac{\left( \omega \left( \left\vert
h\right\vert \right) \right) ^{\frac{\beta }{\eta }}}{\omega ^{\ast }\left(
\left\vert h\right\vert \right) }r^{m}\left( \omega \left( r\right) \right)
^{1-\frac{\beta }{\eta }}
\end{equation*}%
and%
\begin{equation*}
\left\Vert W_{r}\left( f\right) -f\right\Vert _{\omega ^{\ast },~p}=O\left(
1\right) \underset{h\neq 0}{\sup }\frac{\left( \omega \left( \left\vert
h\right\vert \right) \right) ^{\frac{\beta }{\eta }}}{\omega ^{\ast }\left(
\left\vert h\right\vert \right) }r^{m/2}\left( \omega \left( \sqrt{r}\right)
\right) ^{1-\frac{\beta }{\eta }}
\end{equation*}%
as $r\rightarrow 0^{+}$.
\end{cor}

\begin{cor}
Suppose that $0\leq \beta <\alpha \leq 1$, $m\in 
\mathbb{N}
\cup \left\{ 0\right\} $ and $f\in H_{m}^{\alpha ,~p}\left( 1\leq p\leq
\infty \right) $. Then%
\begin{equation*}
\left\Vert P_{r}\left( f\right) -f\right\Vert _{\beta ,~p}=O\left(
r^{m+\alpha -\beta }\right) 
\end{equation*}%
and%
\begin{equation*}
\left\Vert W_{r}\left( f\right) -f\right\Vert _{\beta ,~p}=O\left( r^{\left(
m+\alpha -\beta \right) /2}\right) 
\end{equation*}%
as $r\rightarrow 0^{+}$.
\end{cor}

\section{Proofs of the theorems}

\subsection{Proof of Theorem 1}
For $p=\infty $, (\ref{t1}) was proved in \cite[Theorem 1]{5}. Let $1\leq p<\infty $. Then using (\ref{i2}) and (\ref{i3}) we get%
\begin{equation*}
E_{\lambda }\left( x\right) =\lambda \int\limits_{0}^{\infty }\phi
_{x}\left( t\right) \mathcal{K} \left( \lambda t\right) dt
\end{equation*}%
and%
\begin{equation*}
E_{\lambda }\left( x+h,x\right) =\lambda \int\limits_{0}^{\infty }\left(
\phi _{x+h}\left( t\right) -\phi _{x}\left( t\right) \right) \mathcal{K}
\left( \lambda t\right) dt.
\end{equation*}%
Applying the Fubini inequality \cite{10} we have%
\begin{equation*}
\left\Vert E_{\lambda }\left( \cdot +h,\cdot \right) \right\Vert
_{p}=\lambda \left\{ \int\limits_{%
\mathbb{R}
}\left\vert \int\limits_{0}^{\infty }\left( \phi _{x+h}\left( t\right)
-\phi _{x}\left( t\right) \right) \mathcal{K} \left( \lambda t\right)
dt\right\vert ^{p}dx\right\} ^{\frac{1}{p}}
\end{equation*}%
\begin{equation*}
\leq \lambda \int\limits_{0}^{\infty }\left\vert \mathcal{K} \left( \lambda
t\right) \right\vert \left\{ \int\limits_{%
\mathbb{R}
}\left\vert \phi _{x+h}\left( t\right) -\phi _{x}\left( t\right) \right\vert
^{p}dx\right\} ^{\frac{1}{p}}dt
\end{equation*}%
\begin{equation}
=\lambda \left( \int\limits_{0}^{1/\lambda }+\int\limits_{1/\lambda
}^{1}+\int\limits_{1}^{\infty }\right) \left\vert \mathcal{K} \left(
\lambda t\right) \right\vert \left\{ \int\limits_{%
\mathbb{R}
}\left\vert \phi _{x+h}\left( t\right) -\phi _{x}\left( t\right) \right\vert
^{p}dx\right\} ^{\frac{1}{p}}dt=I_{1}+I_{1}+I_{3}.  \label{p0}
\end{equation}%
It is clear that for $1\leq p<\infty $,%
\begin{equation}
\left\{ \int\limits_{%
\mathbb{R}
}\left\vert \phi _{x+h}\left( t\right) -\phi _{x}\left( t\right) \right\vert
^{p}dx\right\} ^{\frac{1}{p}}\leq 4\omega \left( f,\left\vert h\right\vert
\right) _{p}  \label{p1}
\end{equation}%
and $t\geq 0$%
\begin{equation}
\left\{ \int\limits_{%
\mathbb{R}
}\left\vert \phi _{x+h}\left( t\right) -\phi _{x}\left( t\right) \right\vert
^{p}dx\right\} ^{\frac{1}{p}}\leq 4\omega \left( f,t\right) _{p}.  \label{p2}
\end{equation}%
Then, in view of the property (\ref{i4}) of the kernel $\mathcal{K} $ and
inequalities (\ref{p1}) and (\ref{p2}), we obtain that for $f\in H^{\omega
,~p}$%
\begin{equation*}
I_{1}=\lambda \int\limits_{0}^{1/\lambda }\left\vert \mathcal{K} \left(
\lambda t\right) \right\vert \left\{ \int\limits_{%
\mathbb{R}
}\left\vert \phi _{x+h}\left( t\right) -\phi _{x}\left( t\right) \right\vert
^{p}dx\right\} ^{\frac{1}{p}\frac{\beta }{\eta }}
\end{equation*}%
\begin{equation*}
\cdot \left\{ \int\limits_{%
\mathbb{R}
}\left\vert \phi _{x+h}\left( t\right) -\phi _{x}\left( t\right) \right\vert
^{p}dx\right\} ^{\frac{1}{p}\left( 1-\frac{\beta }{\eta }\right) }dt
\end{equation*}%
\begin{equation*}
=O\left( 1\right) \lambda \left( \omega \left( \left\vert h\right\vert
\right) \right) ^{\frac{\beta }{\eta }}\int\limits_{0}^{1/\lambda
}\left\vert \mathcal{K} \left( \lambda t\right) \right\vert \left( \omega
\left( f,t\right) _{p}\right) ^{1-\frac{\beta }{\eta }}dt
\end{equation*}%
\begin{equation}
=O\left( 1\right) \left( \omega \left( \left\vert h\right\vert \right)
\right) ^{\frac{\beta }{\eta }}\left( \omega \left( \frac{1}{\lambda }%
\right) \right) ^{1-\frac{\beta }{\eta }}.  \label{p3}
\end{equation}

Further, in view of property (\ref{i5}) of the kernel $\mathcal{K} $, by (%
\ref{p1}) and (\ref{p2}) we have that for $f\in H^{\omega ,~p}$%
\begin{equation*}
I_{2}=O\left( 1\right) \lambda \left( \omega \left( \left\vert h\right\vert
\right) \right) ^{\frac{\beta }{\eta }}\int\limits_{1/\lambda
}^{1}\left\vert \mathcal{K} \left( \lambda t\right) \right\vert \left(
\omega \left( f,t\right) _{p}\right) ^{1-\frac{\beta }{\eta }}dt
\end{equation*}%
\begin{equation*}
=O\left( 1\right) \lambda \left( \omega \left( \left\vert h\right\vert
\right) \right) ^{\frac{\beta }{\eta }}\int\limits_{1/\lambda }^{1}\left(
\omega \left( f,t\right) _{p}\right) ^{1-\frac{\beta }{\eta }}\frac{1}{%
\left( \lambda t\right) ^{2}}dt
\end{equation*}%
\begin{equation*}
=O\left( 1\right) \left( \omega \left( \left\vert h\right\vert \right)
\right) ^{\frac{\beta }{\eta }}\frac{1}{\lambda }\int\limits_{1}^{\lambda
}\left( \omega \left( f;\frac{1}{u}\right) _{p}\right) ^{1-\frac{\beta }{%
\eta }}du
\end{equation*}%
\begin{equation}
=O\left( 1\right) \left( \omega \left( \left\vert h\right\vert \right)
\right) ^{\frac{\beta }{\eta }}\frac{1}{\lambda }\int\limits_{1}^{\lambda
}\left( \omega \left( \frac{1}{u}\right) \right) ^{1-\frac{\beta }{\eta }}du.
\label{p4}
\end{equation}

Applying (\ref{i5}), (\ref{p1}) we get for $f\in H^{\omega ,~p}$ that%
\begin{equation*}
I_{3}=\lambda \int\limits_{0}^{1/\lambda }\left\vert \mathcal{K} \left(
\lambda t\right) \right\vert \left\{ \int\limits_{%
\mathbb{R}
}\left\vert \phi _{x+h}\left( t\right) -\phi _{x}\left( t\right) \right\vert
^{p}dx\right\} ^{\frac{1}{p}\frac{\beta }{\eta }}
\end{equation*}%
\begin{equation*}
\cdot \left\{ \int\limits_{%
\mathbb{R}
}\left\vert \phi _{x+h}\left( t\right) -\phi _{x}\left( t\right) \right\vert
^{p}dx\right\} ^{\frac{1}{p}\left( 1-\frac{\beta }{\eta }\right) }dt
\end{equation*}%
\begin{equation*}
=O\left( 1\right) \lambda \left( \omega \left( \left\vert h\right\vert
\right) \right) ^{\frac{\beta }{\eta }}\int\limits_{1}^{\infty }\left\vert
\mathcal{K} \left( \lambda t\right) \right\vert \left( 8\left\Vert
f\right\Vert _{p}\right) ^{1-\frac{\beta }{\eta }}dt
\end{equation*}%
\begin{equation*}
=O\left( 1\right) \lambda \left( \omega \left( \left\vert h\right\vert
\right) \right) ^{\frac{\beta }{\eta }}\int\limits_{1}^{\infty }\frac{1}{%
\left( \lambda t\right) ^{2}}dt=O\left( 1\right) \left( \omega \left(
\left\vert h\right\vert \right) \right) ^{\frac{\beta }{\eta }}\frac{1}{%
\lambda }
\end{equation*}%
\begin{equation*}
=O\left( 1\right) \left( \omega \left( \left\vert h\right\vert \right)
\right) ^{\frac{\beta }{\eta }}\frac{1}{\lambda ^{\frac{\beta }{\eta }}}%
\frac{1}{\lambda ^{1-\frac{\beta }{\eta }}}.
\end{equation*}%
If $\lambda >1$ then $\omega \left( f;1\right) _{p}\leq 2\lambda \omega
\left( f,\frac{1}{\lambda }\right) $. Thus%
\begin{equation*}
I_{3}=O\left( 1\right) \left( \omega \left( \left\vert h\right\vert \right)
\right) ^{\frac{\beta }{\eta }}\left( \omega \left( f,\frac{1}{\lambda }%
\right) _{p}\right) ^{1-\frac{\beta }{\eta }}\frac{1}{\lambda ^{\frac{\beta 
}{\eta }}}
\end{equation*}%
\begin{equation}
=O\left( 1\right) \left( \omega \left( \left\vert h\right\vert \right)
\right) ^{\frac{\beta }{\eta }}\left( \omega \left( \frac{1}{\lambda }%
\right) \right) ^{1-\frac{\beta }{\eta }}.  \label{p5}
\end{equation}%
From relation (\ref{p0}), (\ref{p3}), (\ref{p4}) and (\ref{p5}), we obtain%
\begin{equation*}
\left\Vert E_{\lambda }\left( \cdot +h,\cdot \right) \right\Vert
_{p}=O\left( 1\right) \left( \omega \left( \left\vert h\right\vert \right)
\right) ^{\frac{\beta }{\eta }}\left\{ \left( \omega \left( \frac{1}{\lambda 
}\right) \right) ^{1-\frac{\beta }{\eta }}+\frac{1}{\lambda }%
\int\limits_{1}^{\lambda }\left( \omega \left( \frac{1}{u}\right) \right)
^{1-\frac{\beta }{\eta }}du\right\} .
\end{equation*}%
Since%
\begin{equation*}
\frac{1}{\lambda }\int\limits_{1}^{\lambda }\left( \omega \left( \frac{1}{u}%
\right) \right) ^{1-\frac{\beta }{\eta }}du\geq \frac{1}{\lambda }%
\int\limits_{\frac{\lambda -1}{2}}^{\lambda }\left( \omega \left( \frac{1}{u%
}\right) \right) ^{1-\frac{\beta }{\eta }}du
\end{equation*}%
\begin{equation*}
\geq \left( \omega \left( \frac{1}{\lambda }\right) \right) ^{1-\frac{\beta 
}{\eta }}\frac{\lambda +1}{2\lambda }\geq \frac{1}{2}\left( \omega \left( 
\frac{1}{\lambda }\right) \right) ^{1-\frac{\beta }{\eta }},
\end{equation*}%
therefore%
\begin{equation*}
\left\Vert E_{\lambda }\left( \cdot +h,\cdot \right) \right\Vert
_{p}=O\left( 1\right) \left( \omega \left( \left\vert h\right\vert \right)
\right) ^{\frac{\beta }{\eta }}\left\{ \frac{1}{\lambda }\int\limits_{1}^{%
\lambda }\left( \omega \left( \frac{1}{u}\right) \right) ^{1-\frac{\beta }{%
\eta }}du\right\} .
\end{equation*}%
Hence%
\begin{equation}
\underset{h\neq 0}{\sup }\frac{\left\Vert E_{\lambda }\left( \cdot +h,\cdot
\right) \right\Vert _{p}}{\omega ^{\ast }\left( \left\vert h\right\vert
\right) }=O\left( 1\right) \frac{\left( \omega \left( \left\vert
h\right\vert \right) \right) ^{\frac{\beta }{\eta }}}{\omega ^{\ast }\left(
\left\vert h\right\vert \right) }\left\{ \frac{1}{\lambda }%
\int\limits_{1}^{\lambda }\left( \omega \left( \frac{1}{u}\right) \right)
^{1-\frac{\beta }{\eta }}du\right\} .  \label{p6}
\end{equation}%
We can easily see that%
\begin{equation}
\left\Vert E_{\lambda }\left( \cdot \right) \right\Vert _{p}=O\left(
1\right) \frac{1}{\lambda }\int\limits_{1}^{\lambda }\omega \left( \frac{1}{%
u}\right) du=O\left( 1\right) \frac{1}{\lambda }\int\limits_{1}^{\lambda
}\left( \omega \left( \frac{1}{u}\right) \right) ^{1-\frac{\beta }{\eta }}du.
\label{p7}
\end{equation}

From (\ref{p6}) an (\ref{p7}), we finally obtain%
\begin{equation*}
\left\Vert E_{\lambda }\left( \cdot \right) \right\Vert _{\omega ^{\ast
},~p}=\left\Vert E_{\lambda }\left( \cdot \right) \right\Vert _{p}+\underset{%
h\neq 0}{\sup }\frac{\left\Vert E_{\lambda }\left( \cdot +h,\cdot \right)
\right\Vert _{p}}{\omega ^{\ast }\left( \left\vert h\right\vert \right) }
\end{equation*}%
\begin{equation*}
=O\left( 1\right) \frac{\left( \omega \left( \left\vert h\right\vert \right)
\right) ^{\frac{\beta }{\eta }}}{\omega ^{\ast }\left( \left\vert
h\right\vert \right) }\left\{ \frac{1}{\lambda }\int\limits_{1}^{\lambda
}\left( \omega \left( \frac{1}{u}\right) \right) ^{1-\frac{\beta }{\eta }%
}du\right\} .
\end{equation*}%
This completes the proof of Theorem 1. $\square $

\subsection{Proof of Theorem 2}

Let $p=\infty .$ Then by (\ref{i2}) and (\ref{i3}) we get%
\begin{equation*}
\left\vert E_{\lambda }\left( x+h,x\right) \right\vert \leq \lambda
\int\limits_{0}^{\infty }\left\vert \mathcal{K} \left( \lambda t\right)
\right\vert \left\vert \phi _{x+h}\left( t\right) -\phi _{x}\left( t\right)
\right\vert dt
\end{equation*}%
\begin{equation}
=\lambda \left( \int\limits_{0}^{1/\lambda }+\int\limits_{1/\lambda
}^{\infty }\right) \left\vert \mathcal{K} \left( \lambda t\right)
\right\vert \left\vert \phi _{x+h}\left( t\right) -\phi _{x}\left( t\right)
\right\vert dt=J_{1}+J_{2}.  \label{p8}
\end{equation}%
It is clear that 
\begin{equation*}
\left\vert \phi _{x+h}\left( t\right) -\phi _{x}\left( t\right) \right\vert
\leq 4\omega \left( f;\left\vert h\right\vert \right) _{\infty }
\end{equation*}%
and%
\begin{equation}
\left\vert \phi _{x+h}\left( t\right) -\phi _{x}\left( t\right) \right\vert
\leq 4\omega \left( f;t\right) _{\infty }.  \label{p9}
\end{equation}%
Using this and (\ref{i4}) we have that for $f\in H^{\omega ,~p}$%
\begin{equation*}
J_{1}=\lambda \int\limits_{0}^{1/\lambda }\left\vert \mathcal{K} \left(
\lambda t\right) \right\vert \left( \left\vert \phi _{x+h}\left( t\right)
-\phi _{x}\left( t\right) \right\vert \right) ^{\frac{\beta }{\eta }}\left(
\left\vert \phi _{x+h}\left( t\right) -\phi _{x}\left( t\right) \right\vert
\right) ^{1-\frac{\beta }{\eta }}dt
\end{equation*}%
\begin{equation}
=O\left( 1\right) \lambda \left( \omega \left( \left\vert h\right\vert
\right) \right) ^{\frac{\beta }{\eta }}\int\limits_{0}^{1/\lambda }\left(
\omega \left( f;t\right) _{\infty }\right) ^{1-\frac{\beta }{\eta }%
}dt=O\left( 1\right) \left( \omega \left( \left\vert h\right\vert \right)
\right) ^{\frac{\beta }{\eta }}\left( \omega \left( \frac{1}{\lambda }%
\right) \right) ^{1-\frac{\beta }{\eta }}.  \label{p10}
\end{equation}%
Further, by (\ref{p9}) and (\ref{t2}) we obtain that for $f\in H^{\omega ,~p}$%
\begin{equation*}
J_{2}=\lambda \int\limits_{1/\lambda }^{\infty }\left\vert \mathcal{K}
\left( \lambda t\right) \right\vert \left( \left\vert \phi _{x+h}\left(
t\right) -\phi _{x}\left( t\right) \right\vert \right) ^{\frac{\beta }{\eta }%
}\left( \left\vert \phi _{x+h}\left( t\right) -\phi _{x}\left( t\right)
\right\vert \right) ^{1-\frac{\beta }{\eta }}dt
\end{equation*}%
\begin{equation*}
=O\left( 1\right) \left( \omega \left( \left\vert h\right\vert \right)
\right) ^{\frac{\beta }{\eta }}\lambda \int\limits_{1/\lambda }^{\infty
}\left( \omega \left( f;t\right) _{\infty }\right) ^{1-\frac{\beta }{\eta }%
}\left\vert \mathcal{K} \left( \lambda t\right) \right\vert dt
\end{equation*}%
\begin{equation*}
=O\left( 1\right) \left( \omega \left( \left\vert h\right\vert \right)
\right) ^{\frac{\beta }{\eta }}\lambda \int\limits_{1/\lambda }^{\infty
}\left( \frac{\omega \left( f;t\right) _{\infty }}{t}\right) ^{1-\frac{\beta 
}{\eta }}t^{1-\frac{\beta }{\eta }}\left\vert \mathcal{K} \left( \lambda
t\right) \right\vert dt
\end{equation*}%
\begin{equation*}
=O\left( 1\right) \left( \omega \left( \left\vert h\right\vert \right)
\right) ^{\frac{\beta }{\eta }}\left( \omega \left( f;\frac{1}{\lambda }%
\right) _{\infty }\right) ^{1-\frac{\beta }{\eta }}\lambda ^{2-\frac{\beta }{%
\eta }}\int\limits_{1/\lambda }^{\infty }t^{1-\frac{\beta }{\eta }%
}\left\vert \mathcal{K} \left( \lambda t\right) \right\vert dt
\end{equation*}%
\begin{equation*}
=O\left( 1\right) \left( \omega \left( \left\vert h\right\vert \right)
\right) ^{\frac{\beta }{\eta }}\left( \omega \left( \frac{1}{\lambda }%
\right) \right) ^{1-\frac{\beta }{\eta }}\lambda ^{2}\int\limits_{1/\lambda
}^{\infty }t\left\vert \mathcal{K} \left( \lambda t\right) \right\vert dt
\end{equation*}%
\begin{equation*}
=O\left( 1\right) \left( \omega \left( \left\vert h\right\vert \right)
\right) ^{\frac{\beta }{\eta }}\left( \omega \left( \frac{1}{\lambda }%
\right) \right) ^{1-\frac{\beta }{\eta }}\int\limits_{1}^{\infty
}u\left\vert \mathcal{K} \left( u\right) \right\vert du
\end{equation*}%
\begin{equation}
=O\left( 1\right) \left( \omega \left( \left\vert h\right\vert \right)
\right) ^{\frac{\beta }{\eta }}\left( \omega \left( \frac{1}{\lambda }%
\right) \right) ^{1-\frac{\beta }{\eta }}.  \label{p11}
\end{equation}%
Similarly we can prove that%
\begin{equation*}
\left\Vert E_{\lambda }\left( \cdot \right) \right\Vert _{\infty }=O\left(
1\right) \omega \left( f;\frac{1}{\lambda }\right) _{\infty }=O\left(
1\right) \left( \omega \left( f;\frac{1}{\lambda }\right) _{\infty }\right)
^{1-\frac{\beta }{\eta }}\left( \omega \left( f;\frac{1}{\lambda _{0}}%
\right) _{\infty }\right) ^{\frac{\beta }{\eta }}
\end{equation*}%
\begin{equation}
=O\left( 1\right) \left( \omega \left( \frac{1}{\lambda }\right) \right) ^{1-%
\frac{\beta }{\eta }}.  \label{p12}
\end{equation}%
Hence, by (\ref{p8}), (\ref{p10}), (\ref{p11}) and (\ref{p12})%
\begin{equation*}
\left\Vert E_{\lambda }\left( \cdot \right) \right\Vert _{\omega ^{\ast
},~\infty }=\left\Vert E_{\lambda }\left( \cdot \right) \right\Vert _{\infty
}+\underset{h\neq 0}{\sup }\frac{\left\Vert E_{\lambda }\left( \cdot
+h,\cdot \right) \right\Vert _{\infty }}{\omega ^{\ast }\left( \left\vert
h\right\vert \right) }
\end{equation*}%
\begin{equation*}
=O\left( 1\right) \underset{h\neq 0}{\sup }\frac{\left( \omega \left(
\left\vert h\right\vert \right) \right) ^{\frac{\beta }{\eta }}}{\omega
^{\ast }\left( \left\vert h\right\vert \right) }\left( \omega \left( \frac{1%
}{\lambda }\right) \right) ^{1-\frac{\beta }{\eta }}.
\end{equation*}%
Suppose $1\leq p<\infty $. Then using (\ref{i2}), (\ref{i3}) and the Fubini
inequality \cite{10} we get%
\begin{equation*}
\left\Vert E_{\lambda }\left( \cdot +h,\cdot \right) \right\Vert
_{p}=\lambda \left\{ \int\limits_{%
\mathbb{R}
}\left\vert \int\limits_{0}^{\infty }\left( \phi _{x+h}\left( t\right)
-\phi _{x}\left( t\right) \right) \mathcal{K} \left( \lambda t\right)
dt\right\vert ^{p}dx\right\} ^{\frac{1}{p}}
\end{equation*}%
\begin{equation*}
\leq \lambda \int\limits_{0}^{\infty }\left\vert \mathcal{K} \left( \lambda
t\right) \right\vert \left\{ \int\limits_{%
\mathbb{R}
}\left\vert \phi _{x+h}\left( t\right) -\phi _{x}\left( t\right) \right\vert
^{p}dx\right\} ^{\frac{1}{p}}dt
\end{equation*}%
\begin{equation}
=\lambda \left( \int\limits_{0}^{1/\lambda }+\int\limits_{1/\lambda
}^{\infty }\right) \left\vert \mathcal{K} \left( \lambda t\right)
\right\vert \left\{ \int\limits_{%
\mathbb{R}
}\left\vert \phi _{x+h}\left( t\right) -\phi _{x}\left( t\right) \right\vert
^{p}dx\right\} ^{\frac{1}{p}}dt=S_{1}+S_{2}.  \label{p13}
\end{equation}%
In view of property (\ref{i4}) of the kernel $\mathcal{K} $ and inequalities
(\ref{p1}) and (\ref{p2}), we obtain that for $f\in H^{\omega ,~p}$%
\begin{equation*}
S_{1}=\lambda \int\limits_{0}^{1/\lambda }\left\vert \mathcal{K} \left(
\lambda t\right) \right\vert \left\{ \int\limits_{%
\mathbb{R}
}\left\vert \phi _{x+h}\left( t\right) -\phi _{x}\left( t\right) \right\vert
^{p}dx\right\} ^{\frac{1}{p}\frac{\beta }{\eta }}
\end{equation*}%
\begin{equation*}
\cdot \left\{ \int\limits_{%
\mathbb{R}
}\left\vert \phi _{x+h}\left( t\right) -\phi _{x}\left( t\right) \right\vert
^{p}dx\right\} ^{\frac{1}{p}\left( 1-\frac{\beta }{\eta }\right) }dt
\end{equation*}%
\begin{equation*}
=O\left( 1\right) \lambda \left( \omega \left( \left\vert h\right\vert
\right) \right) ^{\frac{\beta }{\eta }}\int\limits_{0}^{1/\lambda
}\left\vert \mathcal{K} \left( \lambda t\right) \right\vert \left( \omega
\left( f,t\right) _{p}\right) ^{1-\frac{\beta }{\eta }}dt
\end{equation*}%
\begin{equation}
=O\left( 1\right) \left( \omega \left( \left\vert h\right\vert \right)
\right) ^{\frac{\beta }{\eta }}\left( \omega \left( \frac{1}{\lambda }%
\right) \right) ^{1-\frac{\beta }{\eta }}.  \label{p14}
\end{equation}%
Further, in view of property (\ref{t2}) of the kernel $\mathcal{K} $, by (%
\ref{p1}) and (\ref{p2}) we have for $f\in H^{\omega ,~p}$%
\begin{equation*}
S_{2}=O\left( 1\right) \lambda \left( \omega \left( \left\vert h\right\vert
\right) \right) ^{\frac{\beta }{\eta }}\int\limits_{1/\lambda }^{\infty
}\left\vert \mathcal{K} \left( \lambda t\right) \right\vert \left( \omega
\left( f,t\right) _{p}\right) ^{1-\frac{\beta }{\eta }}dt
\end{equation*}%
\begin{equation*}
=O\left( 1\right) \left( \omega \left( \left\vert h\right\vert \right)
\right) ^{\frac{\beta }{\eta }}\lambda \int\limits_{1/\lambda }^{\infty
}\left( \frac{\omega \left( f;t\right) _{p}}{t}\right) ^{1-\frac{\beta }{%
\eta }}t^{1-\frac{\beta }{\eta }}\left\vert \mathcal{K} \left( \lambda
t\right) \right\vert dt
\end{equation*}%
\begin{equation*}
=O\left( 1\right) \left( \omega \left( \left\vert h\right\vert \right)
\right) ^{\frac{\beta }{\eta }}\left( \omega \left( f;\frac{1}{\lambda }%
\right) _{p}\right) ^{1-\frac{\beta }{\eta }}\lambda ^{2-\frac{\beta }{\eta }%
}\int\limits_{1/\lambda }^{\infty }t^{1-\frac{\beta }{\eta }}\left\vert
\mathcal{K} \left( \lambda t\right) \right\vert dt
\end{equation*}%
\begin{equation*}
=O\left( 1\right) \left( \omega \left( \left\vert h\right\vert \right)
\right) ^{\frac{\beta }{\eta }}\left( \omega \left( \frac{1}{\lambda }%
\right) \right) ^{1-\frac{\beta }{\eta }}\lambda ^{2}\int\limits_{1/\lambda
}^{\infty }t\left\vert \mathcal{K} \left( \lambda t\right) \right\vert dt
\end{equation*}%
\begin{equation*}
=O\left( 1\right) \left( \omega \left( \left\vert h\right\vert \right)
\right) ^{\frac{\beta }{\eta }}\left( \omega \left( \frac{1}{\lambda }%
\right) \right) ^{1-\frac{\beta }{\eta }}\int\limits_{1}^{\infty
}u\left\vert \mathcal{K} \left( u\right) \right\vert du
\end{equation*}%
\begin{equation}
=O\left( 1\right) \left( \omega \left( \left\vert h\right\vert \right)
\right) ^{\frac{\beta }{\eta }}\left( \omega \left( \frac{1}{\lambda }%
\right) \right) ^{1-\frac{\beta }{\eta }}.  \label{p15}
\end{equation}%
We can easily see that%
\begin{equation*}
\left\Vert E_{\lambda }\left( \cdot \right) \right\Vert _{p}=O\left(
1\right) \omega \left( f;\frac{1}{\lambda }\right) _{p}=O\left( 1\right)
\left( \omega \left( f;\frac{1}{\lambda _{0}}\right) \right) ^{\frac{\beta }{%
\eta }}\left( \omega \left( f;\frac{1}{\lambda }\right) _{p}\right) ^{1-%
\frac{\beta }{\eta }}
\end{equation*}%
\begin{equation}
=O\left( 1\right) \left( \omega \left( \frac{1}{\lambda }\right) \right) ^{1-%
\frac{\beta }{\eta }}.  \label{p16}
\end{equation}%
Hence, by (\ref{p13}) -(\ref{p16})%
\begin{equation*}
\left\Vert E_{\lambda }\left( \cdot \right) \right\Vert _{\omega ^{\ast
},~p}=\left\Vert E_{\lambda }\left( \cdot \right) \right\Vert _{p}+\underset{%
h\neq 0}{\sup }\frac{\left\Vert E_{\lambda }\left( \cdot +h,\cdot \right)
\right\Vert _{p}}{\omega ^{\ast }\left( \left\vert h\right\vert \right) }
\end{equation*}%
\begin{equation*}
=O\left( 1\right) \underset{h\neq 0}{\sup }\frac{\left( \omega \left(
\left\vert h\right\vert \right) \right) ^{\frac{\beta }{\eta }}}{\omega
^{\ast }\left( \left\vert h\right\vert \right) }\left( \omega \left( \frac{1%
}{\lambda }\right) \right) ^{1-\frac{\beta }{\eta }}.
\end{equation*}

The proof is complete. $\square $

\subsubsection{Proof of Remark 2}

Let $p=\infty $. Then by (\ref{i2}) and (\ref{t3}) we get%
\begin{equation*}
\left\Vert F_{\lambda ,~m}\left( f\right) \right\Vert _{\infty }\leq
\sum\limits_{j=0}^{m}\frac{\left\Vert f^{\left( j\right) }\right\Vert
_{\infty }}{j!}\lambda \int\limits_{%
\mathbb{R}
}\left\vert t^{j}\mathcal{K} \left( \lambda t\right) \right\vert
dt=2\sum\limits_{j=0}^{m}\frac{\left\Vert f^{\left( j\right) }\right\Vert
_{\infty }}{j!}\lambda \int\limits_{0}^{\infty }t^{j}\left\vert \mathcal{K}
\left( \lambda t\right) \right\vert dt
\end{equation*}%
\begin{equation}
=2\sum\limits_{j=0}^{m}\frac{\left\Vert f^{\left( j\right) }\right\Vert
_{\infty }}{j!\lambda ^{j}}\int\limits_{0}^{\infty }u^{j}\left\vert
\mathcal{K} \left( u\right) \right\vert du=O\left( 1\right)
\sum\limits_{j=0}^{m}\frac{\left\Vert f^{\left( j\right) }\right\Vert
_{\infty }}{j!\lambda ^{j}}.  \label{k1}
\end{equation}%
Suppose that $1\leq p<\infty $. Then using the Fubini inequality \cite{10}, (\ref{i2})
and (\ref{t3}) we obtain%
\begin{equation*}
\left\Vert F_{\lambda ,~m}\left( f\right) \right\Vert _{p}=\left\Vert
\sum\limits_{j=0}^{m}\frac{\left( -1\right) ^{j}}{j!}\lambda \int\limits_{%
\mathbb{R}
}f^{\left( j\right) }\left( t+\cdot \right) t^{j}\mathcal{K} \left( \lambda
t\right) dt\right\Vert _{p}
\end{equation*}%
\begin{equation*}
\leq \sum\limits_{j=0}^{m}\frac{\lambda }{j!}\left\{ \int\limits_{%
\mathbb{R}
}\left\vert \int\limits_{%
\mathbb{R}
}f^{\left( j\right) }\left( t+x\right) t^{j}\mathcal{K} \left( \lambda
t\right) dt\right\vert ^{p}dx\right\} ^{\frac{1}{p}}
\end{equation*}%
\begin{equation*}
\leq \sum\limits_{j=0}^{m}\frac{\lambda }{j!}\int\limits_{%
\mathbb{R}
}\left\vert t^{j}\mathcal{K} \left( \lambda t\right) \right\vert \left\{
\int\limits_{%
\mathbb{R}
}\left\vert f^{\left( j\right) }\left( t+x\right) \right\vert ^{p}dx\right\}
^{\frac{1}{p}}
\end{equation*}%
\begin{equation}
=\sum\limits_{j=0}^{m}\frac{\left\Vert f^{\left( j\right) }\right\Vert _{p}%
}{j!}\lambda \int\limits_{0}^{\infty }t^{j}\left\vert \mathcal{K} \left(
\lambda t\right) \right\vert dt=O\left( 1\right) \sum\limits_{j=0}^{m}\frac{%
\left\Vert f^{\left( j\right) }\right\Vert _{p}}{j!\lambda ^{j}}.  \label{k2}
\end{equation}%
Using (\ref{k1}) and (\ref{k2}) we get%
\begin{equation*}
\left\Vert F_{\lambda ,~m}\left( f\right) \right\Vert _{\omega ^{\ast
},~p}=\left\Vert F_{\lambda ,~m}\left( f\right) \right\Vert _{p}+\underset{%
h\neq 0}{\sup }\frac{\left\Vert \Delta _{h}F_{\lambda ,~m}\left( f;\cdot
\right) \right\Vert _{p}}{\omega ^{\ast }\left( \left\vert h\right\vert
\right) }
\end{equation*}%
\begin{equation*}
=\left\Vert F_{\lambda ,~m}\left( f\right) \right\Vert _{p}+\underset{h\neq 0}%
{\sup }\frac{\left\Vert F_{\lambda ,~m}\left( \Delta _{h}f;\cdot \right)
\right\Vert _{p}}{\omega ^{\ast }\left( \left\vert h\right\vert \right) }
\end{equation*}%
\begin{equation*}
=O\left( 1\right) \left\{ \sum\limits_{j=0}^{m}\frac{\left\Vert f^{\left(
j\right) }\right\Vert _{p}}{j!\lambda ^{j}}+\underset{h\neq 0}{\sup }%
\sum\limits_{j=0}^{m}\frac{\left\Vert \Delta _{h}f^{\left( j\right) }\left(
\cdot \right) \right\Vert _{p}}{j!\lambda ^{j}\omega ^{\ast }\left(
\left\vert h\right\vert \right) }\right\} =O\left( 1\right)
\sum\limits_{j=0}^{m}\frac{\left\Vert f^{\left( j\right) }\right\Vert
_{\omega ^{\ast },p}}{j!\lambda ^{j}}.
\end{equation*}

This ends our proof. $\square $

\subsection{Proof of Theorem 3}

We use the following modified Taylor formula for $f\in L_{m}^{p}$ with $m\in 
\mathbb{N}
$:%
\begin{equation*}
f\left( x\right) =\sum\limits_{j=0}^{m}\frac{f^{\left( j\right) }\left(
t\right) }{j!}\left( x-t\right) ^{j}
\end{equation*}%
\begin{equation*}
+\frac{\left( x-t\right) ^{m}}{\left( m-1\right) !}\int\limits_{0}^{1}%
\left( 1-u\right) ^{m-1}\left\{ f^{\left( m\right) }\left( t+u\left(
x-t\right) \right) -f^{\left( m\right) }\left( t\right) \right\} du
\end{equation*}%
for a fixed $t\in 
\mathbb{R}
$ and every $x\in 
\mathbb{R}
$.

By (\ref{i2}) we get%
\begin{equation*}
f\left( x\right) =\lambda \int\limits_{%
\mathbb{R}
}f\left( x\right) \mathcal{K} \left( \lambda \left( t-x\right) \right)
dt=\lambda \int\limits_{%
\mathbb{R}
}\sum\limits_{j=0}^{m}\frac{f^{\left( j\right) }\left( t\right) }{j!}\left(
x-t\right) ^{j}\mathcal{K} \left( \lambda \left( t-x\right) \right) dt
\end{equation*}%
\begin{equation*}
+\lambda \int\limits_{%
\mathbb{R}
}\mathcal{K} \left( \lambda \left( t-x\right) \right) \frac{\left(
x-t\right) ^{m}}{\left( m-1\right) !}
\end{equation*}%
\begin{equation*}
\cdot \left( \int\limits_{0}^{1}\left( 1-u\right) ^{m-1}\left\{ f^{\left(
m\right) }\left( t+u\left( x-t\right) \right) -f^{\left( m\right) }\left(
t\right) \right\} du\right) dt
\end{equation*}%
\begin{equation*}
=F_{\lambda ,m}\left( f;x\right) +\lambda \int\limits_{%
\mathbb{R}
}\mathcal{K} \left( \lambda \left( t-x\right) \right) \frac{\left(
x-t\right) ^{m}}{\left( m-1\right) !}
\end{equation*}%
\begin{equation*}
\cdot \left( \int\limits_{0}^{1}\left( 1-u\right) ^{m-1}\left\{ f^{\left(
m\right) }\left( t+u\left( x-t\right) \right) -f^{\left( m\right) }\left(
t\right) \right\} du\right) dt.
\end{equation*}%
Therefore, by (\ref{i2})%
\begin{equation*}
f\left( x\right) -F_{\lambda ,~m}\left( f;x\right)
\end{equation*}%
\begin{equation*}
=\lambda \int\limits_{%
\mathbb{R}
}\left( \frac{\left( x-t\right) ^{m}}{\left( m-1\right) !}%
\int\limits_{0}^{1}\left( 1-u\right) ^{m-1}\Delta _{u\left( x-t\right)
}f^{\left( m\right) }\left( t\right) du\right) \mathcal{K} \left( \lambda
\left( t-x\right) \right) dt
\end{equation*}%
\begin{equation}
=\lambda \int\limits_{%
\mathbb{R}
}\left( \frac{t^{m}}{\left( m-1\right) !}\int\limits_{0}^{1}\left(
1-u\right) ^{m-1}\Delta _{ut}f^{\left( m\right) }\left( x-t\right) du\right)
\mathcal{K} \left( \lambda t\right) dt.  \label{s00}
\end{equation}%
Set%
\begin{equation*}
E_{\lambda ,~m}\left( x\right) =E_{\lambda ,~m}\left( f;x\right) :=f\left(
x\right) -F_{\lambda ,~m}\left( f;x\right)
\end{equation*}%
and%
\begin{equation*}
E_{\lambda ,~m}\left( x+h,x\right) =E_{\lambda ,~m}\left( f;x\right)
:=E_{\lambda ,~m}\left( x+h\right) -E_{\lambda ,~m}\left( x\right) .
\end{equation*}%
Let $p=\infty $. Then%
\begin{equation*}
\left\vert E_{\lambda ,~m}\left( x+h,x\right) \right\vert \leq \lambda
\int\limits_{%
\mathbb{R}
}\left\vert \mathcal{K} \left( \lambda t\right) \right\vert \frac{\left\vert
t\right\vert ^{m}}{\left( m-1\right) !}
\end{equation*}%
\begin{equation*}
\cdot \left( \int\limits_{0}^{1}\left( 1-u\right) ^{m-1}\left\vert \Delta
_{ut}f^{\left( m\right) }\left( x+h-t\right) -\Delta _{ut}f^{\left( m\right)
}\left( x-t\right) \right\vert du\right) dt.
\end{equation*}%
It is clear that%
\begin{equation*}
\left\vert \Delta _{ut}f^{\left( m\right) }\left( x+h-t\right) -\Delta
_{ut}f^{\left( m\right) }\left( x-t\right) \right\vert \leq 2\omega \left(
f^{\left( m\right) };\left\vert ut\right\vert \right) _{\infty }
\end{equation*}%
and%
\begin{equation*}
\left\vert \Delta _{ut}f^{\left( m\right) }\left( x+h-t\right) -\Delta
_{ut}f^{\left( m\right) }\left( x-t\right) \right\vert \leq 2\omega \left(
f^{\left( m\right) };\left\vert h\right\vert \right) _{\infty }.
\end{equation*}%
This and the properties of the modulus of continuity yields for $f\in
H_{m}^{\omega ,~p}$%
\begin{equation*}
\left\vert E_{\lambda ,~m}\left( x+h,x\right) \right\vert \leq \lambda \left(
\omega \left( f^{\left( m\right) };\left\vert h\right\vert \right) _{\infty
}\right) ^{\frac{\beta }{\eta }}\int\limits_{%
\mathbb{R}
}\frac{\left\vert t\right\vert ^{m}}{\left( m-1\right) !}\left\vert
\mathcal{K} \left( \lambda t\right) \right\vert
\end{equation*}%
\begin{equation*}
\cdot \left( \int\limits_{0}^{1}\left( 1-u\right) ^{m-1}\left( \omega
\left( f^{\left( m\right) };\left\vert ut\right\vert \right) _{\infty
}\right) ^{1-\frac{\beta }{\eta }}du\right) dt
\end{equation*}%
\begin{equation*}
=O\left( 1\right) \lambda \left( \omega \left( \left\vert h\right\vert
\right) \right) ^{\frac{\beta }{\eta }}
\end{equation*}%
\begin{equation*}
\cdot \int\limits_{%
\mathbb{R}
}\left( \frac{\left\vert t\right\vert ^{m}}{\left( m-1\right) !}\left\vert
\mathcal{K} \left( \lambda t\right) \right\vert \left( \omega \left(
f^{\left( m\right) };\left\vert t\right\vert \right) _{\infty }\right) ^{1-%
\frac{\beta }{\eta }}\int\limits_{0}^{1}\left( 1-u\right) ^{m-1}du\right) dt
\end{equation*}%
\begin{equation*}
=O\left( 1\right) \left( \omega \left( \left\vert h\right\vert \right)
\right) ^{\frac{\beta }{\eta }}\left( \omega \left( f^{\left( m\right) };%
\frac{1}{\lambda }\right) _{\infty }\right) ^{1-\frac{\beta }{\eta }}\lambda
\int\limits_{%
\mathbb{R}
}\frac{\left\vert t\right\vert ^{m}}{m!}\left\vert \mathcal{K} \left(
\lambda t\right) \right\vert \left( 1+\lambda \left\vert t\right\vert
\right) ^{1-\frac{\beta }{\eta }}dt.
\end{equation*}%
Using (\ref{i2}) and (\ref{t4})%
\begin{equation*}
=O\left( 1\right) \left( \omega \left( \left\vert h\right\vert \right)
\right) ^{\frac{\beta }{\eta }}\left( \omega \left( f^{\left( m\right) };%
\frac{1}{\lambda }\right) _{\infty }\right) ^{1-\frac{\beta }{\eta }}\frac{%
\lambda }{m!}\int\limits_{0}^{\infty }t^{m}\left( 1+\lambda t\right)
\left\vert \mathcal{K} \left( \lambda t\right) \right\vert dt
\end{equation*}%
\begin{equation*}
=O\left( 1\right) \left( \omega \left( \left\vert h\right\vert \right)
\right) ^{\frac{\beta }{\eta }}\left( \omega \left( \frac{1}{\lambda }%
\right) \right) ^{1-\frac{\beta }{\eta }}\left( \lambda
\int\limits_{0}^{\infty }t^{m}\left\vert \mathcal{K} \left( \lambda
t\right) \right\vert dt+\lambda ^{2}\int\limits_{0}^{\infty
}t^{m+1}\left\vert \mathcal{K} \left( \lambda t\right) \right\vert dt\right)
\end{equation*}%
\begin{equation*}
=O\left( 1\right) \left( \omega \left( \left\vert h\right\vert \right)
\right) ^{\frac{\beta }{\eta }}\left( \omega \left( \frac{1}{\lambda }%
\right) \right) ^{1-\frac{\beta }{\eta }}\frac{1}{\lambda ^{m}}\left(
\int\limits_{0}^{\infty }u^{m}\left\vert \mathcal{K} \left( u\right)
\right\vert du+\int\limits_{0}^{\infty }u^{m+1}\left\vert \mathcal{K}
\left( u\right) \right\vert du\right)
\end{equation*}%
\begin{equation}
=O\left( 1\right) \left( \omega \left( \left\vert h\right\vert \right)
\right) ^{\frac{\beta }{\eta }}\left( \omega \left( \frac{1}{\lambda }%
\right) \right) ^{1-\frac{\beta }{\eta }}\frac{1}{\lambda ^{m}}.  \label{s1}
\end{equation}%
We can easily see that%
\begin{equation}
\left\Vert E_{\lambda ,~m}\left( \cdot \right) \right\Vert _{\infty }=O\left(
1\right) \frac{1}{\lambda ^{m}}\omega \left( f^{\left( m\right) };\frac{1}{%
\lambda }\right) _{\infty }=O\left( 1\right) \frac{1}{\lambda ^{m}}\left(
\omega \left( \frac{1}{\lambda }\right) \right) ^{1-\frac{\beta }{\eta }}.
\label{s2}
\end{equation}%
Hence, by (\ref{s1}) and (\ref{s2})%
\begin{equation*}
\left\Vert E_{\lambda ,~m}\left( \cdot \right) \right\Vert _{\omega ^{\ast
},~\infty }=\left\Vert E_{\lambda ,~m}\left( \cdot \right) \right\Vert
_{\infty }+\underset{h\neq 0}{\sup }\frac{\left\Vert E_{\lambda ,~m}\left(
\cdot +h,\cdot \right) \right\Vert _{\infty }}{\omega ^{\ast }\left(
\left\vert h\right\vert \right) }
\end{equation*}%
\begin{equation*}
=O\left( 1\right) \underset{h\neq 0}{\sup }\frac{\left( \omega \left(
\left\vert h\right\vert \right) \right) ^{\frac{\beta }{\eta }}}{\omega
^{\ast }\left( \left\vert h\right\vert \right) }\frac{1}{\lambda ^{m}}\left(
\omega \left( \frac{1}{\lambda }\right) \right) ^{1-\frac{\beta }{\eta }}.
\end{equation*}

Suppose $1\leq p<\infty $. Using (\ref{s00}) the Fubini inequality \cite{10}, we get%
\begin{equation*}
\left\Vert E_{\lambda ,~m}\left( \cdot +h,\cdot \right) \right\Vert _{p}=%
\frac{\lambda }{\left( m-1\right) !}\left\{ \int\limits_{%
\mathbb{R}
}\left\vert \int\limits_{%
\mathbb{R}
}t^{m}\mathcal{K} \left( \lambda t\right) \right. \right.
\end{equation*}%
\begin{equation*}
\left. \left. \cdot \left( \int\limits_{0}^{1}\left( 1-u\right)
^{m-1}\left( \Delta _{ut}f^{\left( m\right) }\left( x+h-t\right) -\Delta
_{ut}f^{\left( m\right) }\left( x-t\right) \right) du\right) dt\right\vert
^{p}dx\right\} ^{\frac{1}{p}}
\end{equation*}%
\begin{equation*}
\leq \frac{\lambda }{\left( m-1\right) !}\int\limits_{%
\mathbb{R}
}\left\vert t^{m}\mathcal{K} \left( \lambda t\right) \right\vert
\end{equation*}%
\begin{equation*}
\cdot \left\{ \int\limits_{%
\mathbb{R}
}\left\vert \int\limits_{0}^{1}\left( 1-u\right) ^{m-1}\left( \Delta
_{ut}f^{\left( m\right) }\left( x+h-t\right) -\Delta _{ut}f^{\left( m\right)
}\left( x-t\right) \right) du\right\vert ^{p}dx\right\} ^{\frac{1}{p}}dt
\end{equation*}%
\begin{equation*}
\leq \frac{\lambda }{\left( m-1\right) !}\int\limits_{%
\mathbb{R}
}\left\vert t^{m}\mathcal{K} \left( \lambda t\right) \right\vert
\end{equation*}%
\begin{equation*}
\cdot \left( \int\limits_{0}^{1}\left( 1-u\right) ^{m-1}\left\{
\int\limits_{%
\mathbb{R}
}\left\vert \Delta _{ut}f^{\left( m\right) }\left( x+h-t\right) -\Delta
_{ut}f^{\left( m\right) }\left( x-t\right) \right\vert ^{p}dx\right\} ^{%
\frac{1}{p}}du\right) dt.
\end{equation*}

It is clear that $1\leq p<\infty $%
\begin{equation*}
\left\{ \int\limits_{%
\mathbb{R}
}\left\vert \Delta _{ut}f^{\left( m\right) }\left( x+h-t\right) -\Delta
_{ut}f^{\left( m\right) }\left( x-t\right) \right\vert ^{p}dx\right\} ^{%
\frac{1}{p}}\leq 2\omega \left( f^{\left( m\right) };\left\vert
ut\right\vert \right) _{p}
\end{equation*}%
and%
\begin{equation*}
\left\{ \int\limits_{%
\mathbb{R}
}\left\vert \Delta _{ut}f^{\left( m\right) }\left( x+h-t\right) -\Delta
_{ut}f^{\left( m\right) }\left( x-t\right) \right\vert ^{p}dx\right\} ^{%
\frac{1}{p}}\leq 2\omega \left( f^{\left( m\right) };\left\vert h\right\vert
\right) _{p}.
\end{equation*}%
Using this, (\ref{i2}), (\ref{t4}) and the properties of the modulus of
continuity we get that for $f\in H_{m}^{\omega ,~p}$%
\begin{equation*}
\left\Vert E_{\lambda ,~m}\left( \cdot +h,\cdot \right) \right\Vert _{p}\leq 
\frac{\lambda }{\left( m-1\right) !}\left( \omega \left( \left\vert
h\right\vert \right) \right) ^{\frac{\beta }{\eta }}
\end{equation*}%
\begin{equation*}
\cdot \int\limits_{%
\mathbb{R}
}\left( \left\vert t\right\vert ^{m}\left\vert \mathcal{K} \left( \lambda
t\right) \right\vert \int\limits_{0}^{1}\left( 1-u\right) ^{m-1}\left(
\omega \left( f^{\left( m\right) };\left\vert ut\right\vert \right)
_{p}\right) ^{1-\frac{\beta }{\eta }}du\right)
\end{equation*}%
\begin{equation*}
=O\left( 1\right) \frac{\lambda }{\left( m-1\right) !}\left( \omega \left(
\left\vert h\right\vert \right) \right) ^{\frac{\beta }{\eta }}
\end{equation*}%
\begin{equation*}
\cdot \int\limits_{%
\mathbb{R}
}\left( \left\vert t\right\vert ^{m}\left\vert \mathcal{K} \left( \lambda
t\right) \right\vert \left( \omega \left( f^{\left( m\right) };\left\vert
t\right\vert \right) _{p}\right) ^{1-\frac{\beta }{\eta }}\int%
\limits_{0}^{1}\left( 1-u\right) ^{m-1}du\right) dt
\end{equation*}%
\begin{equation*}
=O\left( 1\right) \frac{\lambda }{\left( m-1\right) !}\left( \omega \left(
\left\vert h\right\vert \right) \right) ^{\frac{\beta }{\eta }}
\end{equation*}%
\begin{equation*}
\cdot \int\limits_{%
\mathbb{R}
}\left( \left\vert t\right\vert ^{m}\left\vert \mathcal{K} \left( \lambda
t\right) \right\vert \left( \omega \left( f^{\left( m\right) };\left\vert
t\right\vert \right) _{p}\right) ^{1-\frac{\beta }{\eta }}\int%
\limits_{0}^{1}\left( 1-u\right) ^{m-1}du\right) dt
\end{equation*}%
\begin{equation*}
=\left( \omega \left( \left\vert h\right\vert \right) \right) ^{\frac{\beta 
}{\eta }}\left( \omega \left( f^{\left( m\right) };\frac{1}{\lambda }\right)
_{p}\right) ^{1-\frac{\beta }{\eta }}\lambda \int\limits_{%
\mathbb{R}
}\frac{\left\vert t\right\vert ^{m}}{m!}\left\vert \mathcal{K} \left(
\lambda t\right) \right\vert \left( 1+\lambda \left\vert t\right\vert
\right) ^{1-\frac{\beta }{\eta }}dt.
\end{equation*}%
\begin{equation*}
=O\left( 1\right) \left( \omega \left( \left\vert h\right\vert \right)
\right) ^{\frac{\beta }{\eta }}\left( \omega \left( \frac{1}{\lambda }%
\right) \right) ^{1-\frac{\beta }{\eta }}\frac{1}{\lambda ^{m}}\left(
\int\limits_{0}^{\infty }u^{m}\left\vert \mathcal{K} \left( u\right)
\right\vert du+\int\limits_{0}^{\infty }u^{m+1}\left\vert \mathcal{K}
\left( u\right) \right\vert du\right)
\end{equation*}%
\begin{equation*}
=O\left( 1\right) \left( \omega \left( \left\vert h\right\vert \right)
\right) ^{\frac{\beta }{\eta }}\left( \omega \left( \frac{1}{\lambda }%
\right) \right) ^{1-\frac{\beta }{\eta }}\frac{1}{\lambda ^{m}}.
\end{equation*}%
Similarly we can prove that%
\begin{equation*}
\left\Vert E_{\lambda ,~m}\left( \cdot \right) \right\Vert _{p}=O\left(
1\right) \frac{1}{\lambda ^{m}}\omega \left( f^{\left( m\right) };\frac{1}{%
\lambda }\right) _{p}=O\left( 1\right) \frac{1}{\lambda ^{m}}\left( \omega
\left( \frac{1}{\lambda }\right) \right) ^{1-\frac{\beta }{\eta }}.
\end{equation*}%
Hence for $1\leq p<\infty $%
\begin{equation*}
\left\Vert E_{\lambda ,~m}\left( \cdot \right) \right\Vert _{\omega ^{\ast
},~p}=\left\Vert E_{\lambda ,~m}\left( \cdot \right) \right\Vert _{p}+\underset%
{h\neq 0}{\sup }\frac{\left\Vert E_{\lambda ,~m}\left( \cdot +h,\cdot \right)
\right\Vert _{p}}{\omega ^{\ast }\left( \left\vert h\right\vert \right) }
\end{equation*}%
\begin{equation*}
=O\left( 1\right) \underset{h\neq 0}{\sup }\frac{\left( \omega \left(
\left\vert h\right\vert \right) \right) ^{\frac{\beta }{\eta }}}{\omega
^{\ast }\left( \left\vert h\right\vert \right) }\left( \omega \left( \frac{1%
}{\lambda }\right) \right) ^{1-\frac{\beta }{\eta }}.
\end{equation*}

Thus, the proof is completed. $\square $


\begin{thebibliography}{99}
\bibitem{1} N. I. Akhiezer, Lectures on Approximation Theory , Gostekhizdat,
Moscow (1947) [in Russian]

\bibitem{2} P. L. Butzer and R. J. Nessel, Fourier Analysis and
Approximation, Vol. 1, Academic Press, Basel-Birkh\"{a}user-New York (1971).

\bibitem{3} P. Chandra and R. N. Mohapatra, Degree of approximation of
functions in the H\"{o}lder metric, Acta Math. Hung., 41(1-2) (1983), 67-76.

\bibitem{4} B. Firlej and L. Rempulska, On some singular integrals in H\"{o}%
lder spaces, Math. Nachr., 170 (1994), 93-100.

\bibitem{11} W. H. Hsiang, On the degrees of convergence of Abel and
conjugate Abel sums, Functiones et Approximatio, XII (1982), 83-103.


\bibitem{5} R. A. Lasuriya, Approximation of Functions on the real axis by
Fej\'{e}r type operators in the generalized H\"{o}lder metric, Mathematical
Notes, 81 (4) (2007), 483-488.

\bibitem{6} R. N. Mohapatra and R. S. Rodriguez, On the rate of convergence
of singular integrals for H\"{o}lder continuous functions, Math. Nachr., 149
(1990), 117-124.

\bibitem{7} S. Pr\"{o}ssdorf, Zur Konvergens der Fourierreihen H\"{o}%
lderstiger Funktionen, Math. Nachr., 69 (1975), 7-14.

\bibitem{8} L. Rempulska and Z. Walczak, On Modified Picard and
Gauss-Weierstrass singular integrals, Ukrainian Mathematical Journal, 57
(11) (2005), 1844-1852.

\bibitem{9} T. Singh, The approximation of continuous functions in H\"{o}%
lder metric, Matematichni Vesnik, 43 (3-4) (1991), 111-118.

\bibitem{10} A. Zygmund, Trigonometric series, Vol. I, University Press
(Cambridge, 1959).
\end{thebibliography}
\end{document}